\documentclass{article}
\usepackage{amssymb}
\usepackage{amsfonts}
\usepackage{amsmath}

\begin{document}

\title{On the convergence of bilateral ergodic averages}
\author{Christophe Cuny, Yves Derriennic \\
Universit\'{e} de Bretagne Occidentale}
\date{}
\maketitle

\begin{abstract}
We study the almost sure convergence of bilateral ergodic averages for not
necessarily integrable functions and relate it to the ones of the forward
and backward averages, hence complementing results of Wo\'s and the second
named author. In the case of convergence, using results of Furstenberg on
double recurrence, we prove oscillations of the bilateral ergodic averages
around the limit.
\end{abstract}

\noindent\textit{Keywords:} ergodic averages, almost sure convergence,
double recurrence

\smallskip

\noindent\textit{MSC: } 37A05

\medskip

\textbf{1. Introduction\medskip }

For a stationary random sequence, time running from $-\infty $ to $+\infty ,$
without any assumption of integrability, a question arises: what could be
the link between convergence of symmetric bilateral averages, from time $-n$
to time +$n$, and convergence of unilateral averages, from time $0$ to time $%
n,$ usually considered in the ergodic theorem or the law of large numbers.
The present paper is devoted to this question.

Birkhoff's ergodic theorem asserts that on a measure preserving dynamical
system $(\mathbb{X},$ $\mathcal{X}$, $\mu ,$ $T)$ the averages $A_{n}^{+}f=%
\dfrac{1}{n}\sum_{i=0}^{n-1}f\circ T^{i}$ of $f\in L^{1}$ converge $\mu $
a.e.. Still these averages converge too for many finite measurable functions
which are not integrable. Their characterization is difficult since it
cannot depend only on distributions. Assuming $\mu (\mathbb{X})=1,$ and the
transformation $T$ invertible and ergodic as we shall always do in this
paper, it was proved by Wo\'{s} [7] that backward averages $A_{n}^{-}f=%
\dfrac{1}{n}\sum_{i=0}^{n-1}f\circ T^{-i}$ of a measurable function, not
necessarily integrable but finite, converge a.e. to a finite limit if and
only if forward averages $A_{n}^{+}f$ do and the limits are the same. This
result was given a new approach in [1] where furthermore it was shown that
this equivalence fails if infinite limits are considered: $\lim
A_{n}^{+}f=+\infty $, $\lim \sup A_{n}^{-}f=+\infty ,$ $\lim \inf
A_{n}^{-}f=-\infty $ a.e. may coexist.

Here, in the same context, bilateral averages $B_{n}f=\dfrac{1}{2n+1}%
\sum_{i=-n}^{n}f\circ T^{i}$ are considered for finite measurable functions.
We prove that convergence a.e. of bilateral averages $B_{n}f$ is equivalent
to convergence of forward and backward averages, $A_{n}^{+}f$ and $A_{n}^{-}f
$ together, the limit being the same, finite or not. This result holds even
if we know the convergence of $B_{n}f$ \textit{a priori} only on a subset of 
$\mathbb{X}$ of positive measure. In all other cases $\lim \sup
B_{n}f=+\infty ,$ $\lim \inf B_{n}f=-\infty $ a.e.. In [4] another approach
of the comparison of forward, backward and bilateral averages when they have
finite limits on $\mathbb{X}$, is presented.

In the last part the phenomenon of infinite oscillations around the limit,
well known for unilateral averages ([3], see also [5], [6], [1]) is
established also for bilateral averages; this uses elements of Furstenberg's
multirecurrence theory [2].

Somehow these results show the absence of extra-effect of symmetry or
compensation between future and past for a stationary sequence, even though
the same statistical behavior occurs in both directions.\bigskip

\textbf{2. Preliminary approach\medskip }

A first difficulty appearing in the study of the convergence of symmetric
bilateral averages%
\begin{equation*}
B_{n}f=\dfrac{1}{2n+1}\sum_{i=-n}^{n}f\circ T^{i}=\dfrac{1}{2n+1}%
[f+\sum_{i=1}^{n}(f\circ T^{-i}+f\circ T^{i})],
\end{equation*}

\noindent with $f$ a finite measurable function, is the problem of the
convergence of the residual terms $\dfrac{1}{n}(f\circ T^{-n}+f\circ T^{n})$%
, $\dfrac{1}{n}(f\circ T^{-n}-f\circ T^{n+1})$. The convergence to $0$ of
the first one is a necessary condition for the convergence of the Ces\`{a}ro
averages defining $B_{n}f.$ The convergence to $0$ of the second one is a
necessary condition for the invariance of the limit of $B_{n}f,$ if it
exists, since $B_{n}f-B_{n}f\circ T=\dfrac{1}{2n+1}(f\circ T^{-n}-f\circ
T^{n+1}).$ Thus the problem of the convergence of $\dfrac{1}{n}f\circ T^{n}$
is also present.

Because of the invariance of the measure the convergence in $\mu $-measure
of these residual terms to $0$ is obvious hence for this type of convergence
the invariance of the set of convergence and the invariance of the limit of
averages are easy. The same is true for the $\lim \sup $ in $\mu $-measure
sense (that is always dominated by the $\lim \sup $ a.e.). But for a.e.
convergence these questions are not obvious, in particular the invariance of
the set of convergence of the sequence $B_{n}f$. As will be shown by Lemma 3
below, using Rokhlin towers it is not difficult to build, for any sequence $%
l(n)$ increasing to +$\infty ,$ a non negative measurable function $f$ such
that $\lim \sup \dfrac{1}{l(n)}f\circ T^{n}=+\infty $ a.e.; even with i.i.d.
random variables, averages may converge in probability to $0$ although the
limsup and liminf are infinite.

If \textit{a priori }the convergence a.e. of $B_{n}f$ is known on a non null
subset $U\subset \mathbb{X},$ the convergence holds also in measure on $U$
hence the sequence $B_{n}f$ must converge in measure to a constant on the
whole space $\mathbb{X}$ by ergodicity. Therefore the limit a.e. of $B_{n}f$
on $U$ must be this same constant. Yet, at this point the divergence a.e. on 
$U^{c}$ is not ruled out since the invariance of the set of convergence a.e.
of the sequence $B_{n}f$ remains to be proved. If $\lim B_{n}f$ exists and
is finite a.e. on $U$ then $\lim \dfrac{1}{n}(f\circ T^{-n}+f\circ T^{n})=0$
a.e. on $U;$ if, at this point, we knew that this convergence implies $\lim 
\dfrac{1}{n}f\circ T^{n}=0$ a.e., the invariance of $\lim \dfrac{1}{n}f\circ
T^{n}$ being obvious, the invariance we wish for the set of convergence of $%
B_{n}f$ would follow. But this implication also requires a proof; it will be
given in Part 4.

To overcome this difficulty in the next part, we shall use a direct close
analysis of a typical trajectory of the dynamical system. Such a method was
suggested by Benjamin Weiss to the second named author in a conversation
some years ago; it is somehow similar to Wo\'{s}'s method to prove Theorem 2
in [7].

Before going further we recall briefly here the results of [1] that we shall
need or extend.

We use the usual notation $f^{+}=\max (f,0)$ and $f^{-}=(-f)^{+}.$

\textbf{Proposition 0}. \textit{Let }$f$\textit{\ be a finite measurable
real function defined on }$X$\textit{. If }$F=\sup\limits_{n>0}
\sum_{i=0}^{n-1}f\circ T^{i}<\infty $\textit{\ a.e. the following equality,
called "filling scheme equation" in [1], holds:}%
\begin{equation*}
f=-F^{-}+F^{+}-F^{+}\circ T\text{ \textit{a.e..}}
\end{equation*}%
From this equality the following Statements are deduced (we use the
numbering of Part 4 in [1]):

(b) if $\lim \sup \left\vert A_{n}^{+}f\right\vert <\infty $ a.e. then $\lim
A_{n}^{+}f=\lim A_{n}^{-}f$ a.e.

(c) if $\lim \sup A_{n}^{+}f=c\in \mathbb{R}$ and $\lim \inf
A_{n}^{+}f=-\infty $ a.e. then $\lim \inf A_{n}^{-}f=c$ and $\lim \sup
A_{n}^{-}f=+\infty $ a.e.

(d) if $f\geq 0$, either $\lim \dfrac{1}{n}f\circ T^{n}=0$ a.e. or $\lim
\inf \dfrac{1}{n}f\circ T^{n}=0$ a.e. and $\lim \sup \dfrac{1}{n}f\circ
T^{n}=+\infty $ a.e.

(e) if $\lim \sup A_{n}^{+}f=+\infty ,$\textit{\ }$\lim \inf
A_{n}^{+}f=-\infty $\textit{\ }a.e. it is possible that

\noindent $\lim A_{n}^{-}f=+\infty $\ a.e..

\textbf{Remark. }In [1] Statement (e) was shown only for systems having a
Bernoulli factor, using properties of sequences of i.i.d. random variables.
Proposition 1, below, will show it for any ergodic dynamical system.\bigskip

\textbf{3. Main results\medskip }

To begin with, here is an elementary lemma of additive combinatorics that we
need in the sequel.\medskip

\textbf{Lemma 1.}\textit{\ Let }$a$\textit{\ be a positive integer and }$%
\Delta $\textit{\ a set of integers with}

\noindent $\Delta \subset \left[ a,5a\right) .$ \textit{If }$\mathrm{card}\,
\Delta /4a>7/8$\textit{\ then every integer }$i\in \left[ 5a/2,3a\right) $%
\textit{\ can be written as }$i=j-k$\textit{\ with }$j,k\in \Delta $\textit{%
\ and }$j-2k\geq a.\medskip $

\textbf{Proof. }Consider $i\in \left[ 5a/2,3a\right) $ and the couples $%
(x,x+i)$ with $x\in \left[ a,5a\right) $ and $x+i-2x=i-x\geq a.$ These
inequalities yield at once :

$a\leq x\leq i-a<a+i\leq x+i\leq 2i-a<5a.$

Thus the two intervals $\left[ a,i-a\right] $ and $\left[ a+i,2i-a\right] $
are disjoint, included in $\left[ a,5a\right) ,$ with the same length $%
i-2a+1.$ If $i$ could not be written as claimed by the lemma, $\Delta $
would contain at most one term of each couple $(x,x+i),$ hence $i-2a+1$
integers of the interval $\left[ a,5a\right) $ would be excluded from $%
\Delta ;$ we would get

\noindent $\mathrm{card}\, \Delta \leq 4a-(i-2a+1)<7a/2,$ thus $\mathrm{card}%
\, \Delta /4a<7/8.\diamond \medskip $

Now a key idea of the paper appears in the next lemma.

\textbf{Lemma 2. }\textit{Let }$f$\textit{\ be a finite measurable function.
If the set of points }$x$ \textit{such that }$B_{n}f(x)\leq 0$\textit{\ for
all }$n$\textit{\ large enough has positive measure, i.e. if }$\mu \left[
\cup _{N>0}\cap _{n>N}\left\{ B_{n}f\leq 0\right\} \right] >0$ \textit{then }%
$\sup\limits_{n>0}\sum_{i=0}^{n-1}f\circ T^{i}<\infty $\textit{\ a.e. on the
whole space }$\mathbb{X}.\medskip $

\textbf{Proof. }Let us consider $V=\cap _{n\geq N}\left\{ B_{n}f\leq
0\right\} $ with $N$ large enough to get $\mu \left( V\right) >0.$ By
ergodicity for an integer $p$ large enough, $\mu \left( \cup _{p\geq i\geq
0}T^{-i}V\right) >7/8.$ Let us denote by $\mathcal{W}(x)$ the set of passage
times in $(\cup _{p\geq i\geq 0}T^{-i}V)$ of the orbit of $x\in \mathbb{X}.$

By the ergodic theorem $\lim\limits_{n}\dfrac{1}{n}\mathrm{card}\, (\mathcal{%
W}(x)\cap \left[ 0,n\right) )=\mu \left( \cup _{p\geq i\geq 0}T^{-i}V\right) 
$ for a.e. $x\in $ $\mathbb{X}.$ For such an $x,$ $\lim\limits_{a}\dfrac{1}{%
4a}\mathrm{card}\, (\mathcal{W}(x)\cap \left[ a,5a\right) )=\mu \left( \cup
_{p\geq i\geq 0}T^{-i}V\right) >7/8.$

Now assume that three integers $a,$ $k$ and $j$ are such that $a>N+4p,$ $k$
and $j\in \mathcal{W}(x)$ with $a\leq k<j<5a$ and $j-2k\geq a.$ By
definition of $\mathcal{W}(x)$ there exist two integers $r$ and $s$ between $%
0$ and $p,$ depending on $k$ and $j$ respectively, such that $T^{k+r}x\in V$
and $T^{j+s}x\in V.$ By definition of $V,$ since $k\geq a>N+4p,$ we have for
any integer $\zeta $ such that $\left\vert \zeta \right\vert \leq 2p:$

\begin{equation*}
(\ast )\text{ \ \ \ \ \ \ }\sum_{i=-k-r-\zeta }^{k+r+\zeta }f\circ
T^{i}\left( T^{k+r}x\right) =\sum_{i=-\zeta }^{2k+2r+\zeta }f(T^{i}x)\leq 0.
\end{equation*}
Put $\xi =j-2k+s-2r-\zeta -1;$ since $j-2k\geq a,$ we get $\xi \geq N$ and:%
\begin{equation*}
(\ast \ast )\text{ \ \ \ \ \ }\sum_{i=-\xi }^{\xi }f\circ T^{i}\left(
T^{j+s}x\right) =\sum_{i=-\xi +j+s}^{\xi +j+s}f(T^{i}x)\leq 0.
\end{equation*}%
Since $-\xi +j+s=2k+2r+\zeta +1,$ adding the two previous inequalities $%
(\ast )$ and $(\ast \ast )$ we get $\sum_{i=-\zeta }^{\xi +j+s}f(T^{i}x)\leq
0$ where $\xi +j+s=2(j-k)+2(s-r)-\zeta -1.$ Now we can choose $\zeta
=2(s-r)\in \lbrack -2p,2p]$ and we get%
\begin{equation*}
\sum_{i=-\zeta }^{2(j-k)-1}f(T^{i}x)\leq 0.
\end{equation*}

Then using Lemma 1 for $a$ such that $\dfrac{1}{4a}\mathrm{card}\, (\mathcal{%
W}(x)\cap \left[ a,5a\right) )>7/8$ and $a>N+4p,$ we obtain for every $n\in
\lbrack 5a/2,3a)$ two integers $k$ and $j$ with the properties required
above and $n=j-k$, which leads to the inequalities 
\begin{equation*}
\sum_{i=0}^{2n-1}f(T^{i}x)\leq \sum_{i=-2p}^{2p}\left\vert
f(T^{i}x)\right\vert \text{ and }\sum_{i=0}^{2n}f(T^{i}x)\leq \left\vert
f(x)\right\vert +\sum_{i=-2p}^{2p}\left\vert f(T^{i+1}x)\right\vert .
\end{equation*}

When $a$ increases to $\infty $ the intervals $[5a/2,3a)$ cover a half-line
of integers and the desired result follows at once since $p$ depends only on
the set $V$.$\diamond $\medskip

\textbf{Theorem 1. }\textit{Let }$f$\textit{\ be a finite measurable
function. If }$\lim \sup B_{n}f<+\infty $\textit{\ a.e. on a set of positive
measure then }$\lim A_{n}^{+}f=\lim A_{n}^{-}f=\lim B_{n}f\mathit{=c}$%
\textit{\ a.e. on the whole space }$\mathbb{X}$\textit{\ i.e. these three
sequences converge a.e. to the same constant }$c$ \textit{which may be }$%
-\infty .\medskip $

\textbf{Proof. }Consider a finite constant\textbf{\ }$K$ such that $\mu
\left\{ \lim \sup B_{n}f<K-\epsilon \right\} >0$ with $\epsilon >0.$
Applying Lemma 2 to the function $f-K$ and the transformations $T$ or $%
T^{-1} $ we get $\lim \sup A_{n}^{+}f\leq K$ and $\lim \sup A_{n}^{-}f\leq K$
a.e..

If $\lim \sup A_{n}^{+}f=c>-\infty $ then either $\lim A_{n}^{+}f=c$ or $%
\lim \inf A_{n}^{+}f=-\infty $ a.e. by Statement (b). In the first case $%
\lim A_{n}^{-}f=c$ a.e. for the same reason, hence $\lim B_{n}f=c$ a.e.. The
second case is impossible here since

\noindent $\lim \inf A_{n}^{+}f=-\infty $ would imply $\lim \sup
A_{n}^{-}f=+\infty $ a.e. by Statement (c).

If $\lim A_{n}^{+}f=-\infty $ a.e. two cases are again possible by
Proposition 1 below (improvement of Statement (e)). Either $\lim
A_{n}^{-}f=-\infty $ a.e. or $\lim \sup A_{n}^{-}f=+\infty $ and $\lim \inf
A_{n}^{-}f=-\infty $ a.e.. In the first case obviously $\lim B_{n}f=-\infty $
a.e., and again the second case is here excluded, so the theorem is proved.$%
\diamond $\medskip

Before describing the possible behavior of the three sequences $A_{n}^{+}f,$ 
$A_{n}^{-}f,$ $B_{n}f$ we shall give a complement to Statement (e) of [1]
recalled in Part 2, with a better argument valid for any non trivial ergodic
dynamical system (avoiding the properties of stable probability laws used in
[1]).\medskip

\textbf{Proposition 1.} \textit{On any non atomic ergodic system }$(\mathbb{X%
},\mu ,T)$\textit{\ there exist measurable functions }$f$\textit{\ such that 
}$\lim \sup A_{n}^{+}f=+\infty ,$\textit{\ }$\lim \inf A_{n}^{+}f=-\infty $%
\textit{\ a.e. but }$\lim A_{n}^{-}f=+\infty $\textit{\ a.e..}

For the proof we need the following lemma.\medskip

\textbf{Lemma 3.} \textit{On any non atomic ergodic system }$(\mathbb{X},\mu
,T)$ \textit{and for every increasing sequence of integers }$l(n)$\textit{\
there exist measurable functions }$v\geq 0$ \textit{such that }$\lim \sup 
\dfrac{1}{l(n)}v\circ T^{n}=+\infty $\textit{\ a.e..}

\textbf{Proof. }By Rokhlin's lemma (see [6] p.48), for every integer $n$ we
can build a tower of height $n^{2}$ with basis $U_{n}$ and complementary
part $V_{n}$ such that $\sum_{n}\mu (V_{n})<\infty :$ the sets $T^{i}U_{n}$
are pairwise disjoint for $0\leq i\leq n^{2}$ and $V_{n}=(\cup _{0\leq i\leq
n^{2}}U_{n})^{c}.$ Put $g_{n}=nl(n^{2})$ on the roof of the $n^{th}$ tower $%
T^{n^{2}}U_{n}$ and $g_{n}=0$ elsewhere. Put $v=\sum_{n}g_{n};$ it is finite
a.e. by Borel-Cantelli lemma since $\mu (U_{n})\leq n^{-2}.$ Again by
Borel-Cantelli, a.e. $x\in \mathbb{X}$ belongs to all the towers except a
finite number of them since $\sum_{n}\mu (V_{n})<\infty $. If $x$ belongs to
the $n^{th}$ tower there is an integer $i\leq n^{2}$ such that $T^{i}x\in
T^{n^{2}}U_{n}$; thus $\max_{j\leq n^{2}}v(T^{j}x)\geq \max_{j\leq
n^{2}}g_{n}(T^{j}x)\geq nl(n^{2})$. Therefore $l(n^{2})^{-1}\max_{j\leq
n^{2}}v(T^{j}x)\geq n$ a.e. for $n$ large enough, thus $\lim \sup \dfrac{1}{%
l(n)}\max_{j\leq n}v(T^{j}x)=+\infty $ a.e.. Since the sequence $l(n)$ is
increasing, this implies at once the desired result.$\diamond $\medskip

\textbf{Proof of Proposition 1.} An example of a function $f$ having the
desired properties will be given by $f=u+v-v\circ T$ where $u$ and $v$ are
non negative functions such that $\int ud\mu =+\infty $ and $\lim \sup
\Big( \dfrac{1}{n}v\circ T^{n}-A_{n}^{+}u\Big) =+\infty $ a.e.. Indeed, in this case we get $\lim A_{n}^{-}f=\lim (A_{n}^{-}u+\dfrac{1}{n}%
v\circ T^{n})=+\infty \mathit{\ }$a.e.. but $\lim \sup A_{n}^{+}f=\lim \sup
(A_{n}^{+}u-\dfrac{1}{n}v\circ T^{n})=+\infty $ a.e. (recall that since $v\ge 0$ is finite   a.e.  $\liminf \dfrac{1}{n}v\circ T^{n}=0$ a.e.)  and 
$\lim \inf
A_{n}^{+}f=-\lim \sup \Big( \dfrac{1}{n}v\circ T^{n}-A_{n}^{+}u\Big)
=-\infty  \quad \mbox{ a.e..}
$

Let us build two functions $u$ and $v$ with these properties.  By Lemma 3 there is a function $v\geq 0$ such that $\lim \sup n^{-2}(v\circ
T^{n})=+\infty $ a.e..

Let $u$ be a non-negative function such that $u\notin L^{1}$ but $\sqrt{u}%
\in L^{1}.$ By Birkhoff's ergodic theorem we know $\lim \dfrac{1}{n}%
\sum_{0}^{n-1}\sqrt{u}\circ T^{i}=l$ finite a.e.. Now for a non negative
sequence $x_{i},$ if $\dfrac{1}{n}\sum_{0}^{n-1}x_{i}\longrightarrow l,$
then $\dfrac{x_{n}}{n}\longrightarrow 0$ and%
\begin{equation*}
\dfrac{1}{n^{2}}\sum_{0}^{n-1}x_{i}^{2}\leq \dfrac{\max\limits_{i\leq n}x_{i}%
}{n}\left( \dfrac{1}{n}\sum_{0}^{n-1}x_{i}\right) \longrightarrow 0.
\end{equation*}%
Thus $\lim \dfrac{1}{n^{2}}\sum_{0}^{n-1}u\circ T^{i}=0$ a.e..

For these two functions $u$ and $v$, $\lim \sup \dfrac{1}{n}\left( \dfrac{1}{%
n}v\circ T^{n}-A_{n}^{+}u\right) =+\infty $ a.e. and \textit{a fortiori} $%
\lim \sup \left( \dfrac{1}{n}v\circ T^{n}-A_{n}^{+}u\right) =+\infty $ a.e..
So $f=u+v-v\circ T$ satisfies our claim.$\diamond \medskip $

Now we can describe synthetically the possible asymptotic behaviors of the
sequences $A_{n}^{+}f,$ $A_{n}^{-}f,$ $B_{n}f,$ considered together.\medskip

\textbf{Theorem 2}.\textit{\ On any ergodic dynamical system }$(\mathbb{X}%
,\mu ,T),$\textit{\ for a finite measurable function }$f$\textit{\ there are
only two typical situations: }

\textit{1) }$\lim A_{n}^{+}f=\lim A_{n}^{-}f=\lim B_{n}f=c$\textit{\ \ a.e.,
i.e. the three sequences converge a.e. to the same constant which may be }$%
\pm \infty .$

2)$\left\{ 
\begin{array}{l}
with\mathit{\ }c\mathit{\ any\ real\ value\ or}\text{ }+\infty \\ 
\lim \sup A_{n}^{+}f=c,\text{ }\lim \inf A_{n}^{+}f=-\infty \\ 
\lim \sup A_{n}^{-}f=+\infty ,\text{ }\lim \inf A_{n}^{-}f=c \\ 
\lim \sup B_{n}f=+\infty ,\text{ }\lim \inf B_{n}f=-\infty \ \text{\textit{%
a.e.}}%
\end{array}%
\right. $

\textit{(of course the situations derived from 2) by symmetry or sign
changing can also occur).} \thinspace

In particular we see that the sequence of bilateral averages $B_{n}f$
converges a.e. to a finite or infinite limit if and only if the two
sequences of unilateral averages $A_{n}^{+}f$ and $A_{n}^{-}f$ converge to
the same limit. Moreover when the sequence $B_{n}f$ does not converge it
cannot have a finite $\lim \sup $ or $\lim \inf ,$ contrarily to $A_{n}^{+}f$
or $A_{n}^{-}f.$ If the system is non atomic taking two well chosen
nonnegative functions $u$ and $v$ we have case 2) for $f=u-u\circ T+v\circ
T-v$ with $c=+\infty $ and the convergence in measure to $0$ of the three
sequences.

The proof of Theorem 2 is an easy application of Theorem 1 together with
results of [1] recalled in Part 2.\bigskip 

\textbf{4. Complementary results\medskip }

It is natural to consider also asymmetric bilateral averages

$B_{n}^{q}f=\dfrac{1}{2n+1+q}\sum_{i=-n}^{n+q}f\circ T^{i}$ with $q$ a fixed
integer. Recall that \textit{a} \textit{priori }the convergence of $\dfrac{1%
}{n}f\circ T^{n}$ is unknown. Lemma 2\textit{\ }can be reformulated for the
asymmetric averages $B_{n}^{q}f$ with an entirely similar proof. Thus we get
the following result parallel to Theorem 1:\medskip

\textbf{Theorem 3}. \textit{Let }$f$\textit{\ be a finite measurable
function. If }$\lim \sup B_{n}^{q}f<+\infty $\textit{\ a.e. on a set of
positive measure then}

$\lim A_{n}^{+}f=\lim A_{n}^{-}f=\lim B_{n}f\mathit{=\lim B_{n}^{q}f=c}$%
\textit{\ a.e. on the whole space }$\mathbb{X}$\textit{\ with }$\mathit{c}%
\in \mathbb{R}$\textit{\ or }$c=-\infty .$ \textit{Convergences a.e. of }$%
B_{n}f$\textit{\ and }$B_{n}^{q}f$\textit{\ are equivalent.\medskip }

Let us consider now the "residual" terms.\textit{\ A priori }the invariance
of the set of convergence a.e. of $B_{n}f$ was a problem since the
convergence to 0 of $\dfrac{1}{n}(f\circ T^{-n}-f\circ T^{n+1})$ was unknown
even on a part where the convergence of $B_{n}f$ could have been
established. Now we can prove the following propositions.\medskip

\textbf{Proposition 2.}\textit{\ For a finite measurable function }$f$%
\textit{\ and a fixed integer }$q,$\textit{\ if }$\lim \sup \dfrac{1}{n}%
\left\vert f\circ T^{-n}-f\circ T^{n+q}\right\vert <+\infty $\textit{\ a.e.
on a set of} \textit{positive measure} \textit{then }$\lim \dfrac{1}{n}%
f\circ T^{n}=0$\textit{\ a.e. on }$\mathbb{X}$\textit{.}

\textbf{Proof.} Consider $B_{n}^{q-1}(f-f\circ T)=\dfrac{1}{2n+q}(f\circ
T^{-n}-f\circ T^{n+q}).$ Then the result follows from Theorem 3 since $%
A_{n}^{+}(f-f\circ T)=\dfrac{1}{n}(f-f\circ T^{n})$.$\diamond \medskip $

It appears less direct to get the same conclusion starting with

\noindent $\dfrac{1}{n}(f\circ T^{-n}+f\circ T^{n})$. We do not know how to
deduce directly the next proposition from the results of Part 3. We shall
give a proof similar to the one of Lemma 2.\medskip

\textbf{Proposition 3.}\textit{\ Let }$f$ \textit{be} \textit{a finite
measurable function.}

\noindent \textit{If }$\lim \sup \dfrac{1}{n}\left\vert f\circ T^{-n}+f\circ
T^{n}\right\vert <+\infty $\textit{\ a.e. on a set of positive measure} 
\textit{then }$\lim \dfrac{1}{n}f\circ T^{n}=0$\textit{\ a.e}. \textit{on} $%
\mathbb{X}$\textit{.}

First we need a lemma similar to Lemma 1 and even simpler, so we skip its
proof.\medskip

\textbf{Lemma 4. }\textit{Let }$a$\textit{\ be a positive integer and }$%
\Delta $\textit{\ a set of integers with}

\noindent $\Delta \subset \left[ 2a,3a\right) .$ \textit{If }$\mathrm{card}%
\, $ $\Delta /a>3/4$\textit{\ then every integer }$i\in \left[
a/4,a/2\right) $\textit{\ can be written as }$i=j-k$\textit{\ with }$j,k\in
\Delta .\medskip $

\textbf{Proof of Proposition 3. }With $M$ and $N$ large enough put

\noindent $V=\cap _{n>N}\left\{ \dfrac{1}{n}\left\vert f\circ T^{-n}+f\circ
T^{n}\right\vert <M\right\} $ and take $p$ such that

\noindent $\mu \left( \cup _{0\leq i\leq p}T^{-i}V\right) >3/4.$ The set of
passage times $\mathcal{W}(x)$ of the orbit of $x$ in $\left( \cup _{0\leq
i\leq p}T^{-i}V\right) \mathbb{\ }$ has an asymptotic density greater than $%
3/4.$

Assume that three integers $a,k$ and $j$ satisfy: $a>N+p,$ $k$ and $j\in 
\mathcal{W}(x)$ with $2a\leq k<j<3a.$ By definition of $\mathcal{W}(x)$
there exist non negative integers $r,s\leq p$ such that $T^{k+r}x\in V$ and $%
T^{j+s}x\in V.$ For integers $\zeta $ and $\xi >N$ we have%
\begin{equation*}
\left\vert f(T^{k+r-\zeta }x)+f(T^{k+r+\zeta }x)\right\vert <\zeta M
\end{equation*}

and%
\begin{equation*}
\left\vert f(T^{j+s-\xi }x)+f(T^{^{j+s+\xi }}x)\right\vert <\xi M.
\end{equation*}

\noindent The solutions $\zeta $ and $\xi $ of $k+r+\zeta =j+s+\xi $ and $%
j+s-\xi =2(j-k),$ are $\xi =2k-j+s>a>N$ and $\zeta =\xi +j+s-k-r>a-r>N.$
Hence the difference of the terms appearing in the two previous inequalities
yields:%
\begin{equation*}
\left\vert f(T^{2(r-s)}x)-f(T^{2(j-k)}x)\right\vert <(\zeta +\xi )M<(6a+3p)M
\end{equation*}

\noindent since $k+r-\zeta =2(r-s)$ and $\zeta +\xi <6a+3p.$

By Lemma 4, since density of $\mathcal{W}(x)>3/4,$ for all $a$ large enough
every $n\in \lbrack a/4,$ $a/2)$ can be written as $n=j-k$ with $k,$ $j\in 
\mathcal{W}(x)$ and $2a<k<j<3a.$ Thus for such an $n$ we get%
\begin{equation*}
\dfrac{1}{n}\left\vert f(T^{2n}x)\right\vert <\dfrac{4}{a}\left[
(6a+3p)M+\sum_{i=-2p}^{2p}\left\vert f(T^{i}x)\right\vert \right] .
\end{equation*}

When $a$ increases to $\infty $ the intervals $[a/4,a/2)$ cover a half-line
of integers and we obtain%
\begin{equation*}
\lim \sup \dfrac{1}{n}\left\vert f\circ T^{2n}\right\vert <24M<\infty \text{
a.e.}
\end{equation*}%
since $p$ depends only on the set $V$. By Statement (d), the desired result
follows.$\diamond \medskip $

To conclude this part we extend to symmetric bilateral sums the (well-known)
result that $\sup\limits_{n\geq 0}\left\vert \sum_{i=0}^{n}f\circ
T^{i}\right\vert <+\infty $\textit{\ }a.e. implies that $f$ is a "bounded
coboundary", that is $f=g-g\circ T$ with $g$ a bounded measurable function.
Since a precise reference seems difficult to give, we mention here that this
result follows easily from Proposition 0: from $\sup\limits_{n\geq
0}\left\vert \sum_{i=0}^{n}f\circ T^{i}\right\vert <+\infty $\textit{\ }a.e.
we deduce $f=-F^{-}+F^{+}-F^{+}\circ T$ a.e.; then $\lim A_{n}^{+}f=0 $ a.e.
yields $F^{-}=0$ and $F^{+}$ must be bounded since, otherwise, $%
\sup\limits_{n\geq 0}\left\vert F^{+}-F^{+}\circ T^{n}\right\vert =+\infty $%
\textit{\ }a.e. by ergodicity.\medskip

\textbf{Proposition 4}. \textit{Let }$f$ \textit{be} \textit{a finite
measurable function. If }

\noindent $\sup\limits_{n\geq 0}\left\vert \sum_{i=-n}^{n}f\circ
T^{i}\right\vert <+\infty $\textit{\ a.e. on a set of positive measure then
there exists a bounded measurable function }$g$ \textit{such that} $%
f=g-g\circ T$\textit{.}

\textbf{Proof}. Let $K$ be large enough for $\mu \left\{ \sup\limits_{n\geq
0}\left\vert \sum_{i=-n}^{n}f\circ T^{i}\right\vert <K\right\} >0.$
Following the same argument as in Lemma 2 we get for a.e. $x\in \mathbb{X}$
and all $n$ large enough $\left\vert \sum_{i=0}^{2n-1}f(T^{i}x)\right\vert
<2K+\sum_{i=-2p}^{2p}\left\vert f(T^{i}x)\right\vert ,$ where $p$ is
independent of $x$ and $n.$ Hence $\sup\limits_{n\geq 0}\left\vert
\sum_{i=0}^{n}f\circ T^{i}\right\vert <+\infty $\textit{\ }a.e\textit{. }on $%
\mathbb{X}$. Then the conclusion follows from the  result we just
recalled above.$\diamond \bigskip $

\textbf{5. Infinite oscillations around the limit.\medskip }

For unilateral ergodic averages $A_{n}^{+}f=\dfrac{1}{n} \sum_{i=0}^{n-1}f%
\circ T^{i}$ it is well known that infinite oscillations around the limit
must occur; this was first established for $L^{1}-$functions ([3], see also
[5] \S 1.6.3 or [6]). In [1] it is proved for any measurable function for
which the averages converge: the problem is first reduced to the case where $%
f$ is a coboundary, that is of the form $f=g-g\circ T$ with $g$ measurable,
and then the result is an easy by-product of Poincar\'{e}'s recurrence
theorem. Here, for bilateral ergodic averages $B_{n}f$, we shall follow the
same method but the conclusion will require some results of double
recurrence due to Furstenberg [2].\medskip

\textbf{Theorem 4.} \textit{Let }$f$ \textit{be} \textit{a finite measurable
function. If }$\lim B_{n}f=c$\textit{\ a.e. with }$c$\textit{\ finite, then
the sequence }$B_{n}f$\textit{\ oscillates infinitely often around }$c$%
\textit{\ a.e. (in the wide sense; that is the difference }$B_{n}f-c$\textit{%
\ cannot be ultimately strictly positive or ultimately strictly negative on
a set of positive measure).}

\textbf{Proof}. Put $c=0.$ By Theorem 1, $\lim A_{n}^{+}f=0$ a.e..

\textit{Ad absurdum, }suppose that for a.e. $x$ in a set of positive measure$%
,$ there exists $N$ such that $B_{n}f(x)<0$ for all $n>N.$ Then by Lemma 2,
we get $F=\sup\limits_{n>0}\sum_{i=0}^{n-1}f\circ T^{i}<\infty $\ a.e. on $%
\mathbb{X}$ and by Proposition 0, $f=-F^{-}+F^{+}-F^{+}\circ T$ a.e.. Since $%
0=\lim A_{n}^{+}f=-\int F^{-}d\mu $ we obtain $F^{-}=0$ a.e.. Therefore $f$
must be a coboundary: $f=F^{+}-F^{+}\circ T$ and $B_{n}f=\dfrac{1}{2n+1}%
[F^{+}\circ T^{-n}-F^{+}\circ T^{n+1}].$

But in this case it appears a contradiction with the following theorem, that
we shall prove next, which completes the proof of Theorem 4.$\diamond
\medskip $

This last theorem does not depend on the previous results, and might be of
independent interest.\medskip

\textbf{Theorem 5. }\textit{Let }$f$\textit{\ be a measurable function on an
ergodic invertible dynamical system }$(\mathbb{X},$ $\mu ,$ $T).$\textit{\
The set of points }$x$\textit{\ for which the strict inequality}

\noindent $f(T^{-n}x)<f(T^{n}x)$\textit{\ holds for all }$n$\textit{\ large
enough, is negligible, that is}%
\begin{equation*}
\mu \lbrack \cup _{N>0}\cap _{n>N}\left\{ f\circ T^{-n}<f\circ T^{n}\right\}
]=0.
\end{equation*}

\noindent \textit{The same statement holds for inequalities }$%
f(T^{-n}x)<f(T^{n+q}x)$ \textit{with q fixed.\medskip }

Before giving the proof we recall the elements of Furstenberg's
multirecurrence theory for $T$ and $T^{2}$ that we shall need $([2]$, see
the first three parts).

On an ergodic dynamical system $(\mathbb{X},$ $\mu ,$ $T)$ the averages $%
\dfrac{1}{n}\sum_{i=0}^{n-1}(u\circ T^{i})(v\circ T^{2i})$ converge in $%
L^{2}(\mu )$ for every $u$ and $v$ bounded; to represent the limit the
notion of maximal Kronecker factor is introduced. A factor of the ergodic
dynamical system $(\mathbb{X},$ $\mu ,$ $T)$ is a Kronecker factor if it is
a system $(G,$ $m$, $\alpha )$ where $G$ is a compact abelian "monothetic"
group, $m$ its Haar measure and $\alpha \in G$ acts on $G$ by translation $%
\tau _{\alpha }:x\rightarrow x+\alpha $, the sequence $(n\alpha )_{n>0}$
being dense in $G.$ For the maximal Kronecker factor there exists a positive
linear operator $\pi $ of $L^{2}(\mathbb{X},\mu )$ onto $L^{2}(G,m)$ with $%
\pi (f\circ T)=\pi (f)\circ \tau _{\alpha },$ preserving the integral, such
that for every bounded $u,$ $v$ and $w$ the following identity, called
Furstenberg's identity, holds:%
\begin{equation*}
\lim_{n}\dfrac{1}{n}\sum_{i=0}^{n-1}\int_{\mathbb{X}%
}u(x)v(T^{i}x)w(T^{2i}x)d\mu (x)
\end{equation*}%
\begin{equation*}
=\int \int_{G\times G}\pi u(z)\pi v(z+z^{\prime })\pi w(z+2z^{\prime
})dm(z)dm(z^{\prime }).
\end{equation*}

\noindent Here $2z=z+z$ defines a continuous endomorphism of the group $G$
that we'll denote $\theta $. The unique ergodicity of the Kronecker system
is essential (see [2] part 3, especially Lemma 3.4).\medskip

\textbf{Proof of Theorem 5. }If $\mu \left\{ f=ess\sup f\right\} >0$ the
result is easy. Indeed by ergodicity $T^{-n}x\in \left\{ f=ess\sup f\right\} 
$ for infinitely many $n$ and a.e. $x$; for them $f(T^{-n}x)\geq f(T^{n}x).$
If $\mu \left\{ f=ess\inf f\right\} >0$ the dual argument applies.

Therefore we may assume $\mu \left\{ f=ess\inf f\right\} =0$ or $f$
unbounded below. In both cases there exist $a\in \mathbb{R}$ for which the
sets $\left\{ f\leq a\right\} $ have positive and arbitrarily small $\mu $%
-measure.

Now let us take a fixed value $a\in \mathbb{R}$ with $0<\mu \left\{ f\leq
a\right\} <1$ (for the conclusion at the end of the proof we shall have to
let $\mu \left\{ f\leq a\right\} \rightarrow 0).$

Let $E$ be the set considered in the statement :%
\begin{equation*}
E=\cup _{N>0}\cap _{n>N}\left\{ f\circ T^{-n}<f\circ T^{n}\right\} .
\end{equation*}

\noindent Clearly if $x\in E,$ for all $n$ large enough $T^{-n}x\in \left\{
f>a\right\} $ implies $T^{n}x\in \left\{ f>a\right\} $. Thus $E\subset \cup
_{N>0}E_{N}^{a}$ where 
\begin{equation*}
E_{N}^{a}=\cap _{n>N}(T^{n}\left\{ f\leq a\right\} \cup T^{-n}\left\{
f>a\right\} ).
\end{equation*}%
For every $n>N,$ using the invariance of the measure we get :%
\begin{equation*}
\mu \left[ E_{N}^{a}\cap T^{n}\left\{ f>a\right\} \cap T^{-n}\left\{ f\leq
a\right\} \right]
\end{equation*}%
\begin{equation*}
=\mu \left[ T^{-n}E_{N}^{a}\cap \left\{ f>a\right\} \cap T^{-2n}\left\{
f\leq a\right\} \right] =0.
\end{equation*}

Now we shall write Furstenberg's identity with $u=1-1_{\left\{ f\leq
a\right\} },$ $v=1_{E_{N}^{a}},$ $w=1_{\left\{ f\leq a\right\} }.$ The two
functions $\varphi _{N}^{a}=\pi v=\pi 1_{E_{N}^{a}}$ and $\psi ^{a}=\pi
w=\pi 1_{\left\{ f\leq a\right\} }$ belong to $L^{2}(G,m)$ and satisfy $%
0\leq \varphi _{N}^{a}\leq 1$ $m$-a.e., $\int_{G}\varphi _{N}^{a}dm=\mu
(E_{N}^{a})$, and $0\leq \psi ^{a}\leq 1$ $m$-a.e., $\int_{G}\psi ^{a}dm=\mu
\left\{ f\leq a\right\} $. The identity together with the preceding equality
yield:%
\begin{equation*}
0=\lim_{n}\dfrac{1}{n}\sum_{i=0}^{n-1}\mu \lbrack T^{-i}E_{N}^{a}\cap
\left\{ f>a\right\} \cap T^{-2i}\left\{ f\leq a\right\} ]
\end{equation*}%
\begin{equation*}
=\int \int_{G\times G}(1-\psi ^{a}(z))\varphi _{N}^{a}(z+z^{\prime })\psi
^{a}(z+2z^{\prime })dm(z)dm(z^{\prime })
\end{equation*}

\noindent for every integer $N.$

Now put $\varphi =\pi 1_{E}$; it satisfies $\varphi \leq \pi
1_{E\smallsetminus E_{N}^{a}}+\varphi _{N}^{a}.$ Since $E\subset \cup
_{N>0}E_{N}^{a},$ the sequence $(E_{N}^{a})_{N}$ being increasing, we have

\noindent $\lim \int \pi 1_{E\smallsetminus E_{N}^{a}}dm=\lim \mu $ $%
(E\smallsetminus E_{N}^{a})=0$, therefore we get%
\begin{equation*}
\int \int_{G\times G}(1-\psi ^{a}(z))\varphi (z+z^{\prime })\psi
^{a}(z+2z^{\prime })dm(z)dm(z^{\prime })=0.
\end{equation*}

Then Haar measure properties yield :%
\begin{equation*}
\int_{G}\varphi (y)\int_{G}\psi ^{a}(z)(1-\psi ^{a}(2y-z))dm(z)dm(y)=0.
\end{equation*}

Now we shall argue by contradiction. Suppose $\mu (E)>0,$ which implies $%
m\left\{ \varphi >0\right\} >0.$ Since $0\leq \varphi $ and $0\leq \psi
^{a}\leq 1$, for $m$-a.e. $y\in \left\{ \varphi >0\right\} $ we get $\psi
^{a}(z)(1-\psi ^{a}(2y-z))=0$ for $m$-a.e. $z\in G.$ Hence $\psi ^{a}(z)>0$
implies $\psi ^{a}(2y-z)=1,$ that is $\left\{ \psi ^{a}>0\right\} \subset
(-A^{a})+2y$ \ where we denote $\left\{ \psi ^{a}=1\right\} =A^{a}.$ Thus $%
m(\left\{ \psi ^{a}>0\right\} )=m(A^{a})$ and $\psi ^{a}=1_{A^{a}}$ $m$%
-a.e.; moreover

\noindent $A^{a}=(-A^{a})+2y.$ In other words $A^{a}$ is invariant by the
map $z\rightarrow 2y-z$ for $m$-a.e. $y\in \left\{ \varphi >0\right\} .$
Taking differences $A^{a}$ is invariant by translations by elements of the
set $(2\left\{ \varphi >0\right\} )-(2\left\{ \varphi >0\right\} )$ and also
of the generated subgroup . This subgroup is $2H=\theta H$ if we denote by $%
H $ the subgroup generated by the set $\left\{ \varphi >0\right\} -\left\{
\varphi >0\right\} .$ By our supposition $m\left\{ \varphi >0\right\} >0,$
the set $\left\{ \varphi >0\right\} -\left\{ \varphi >0\right\} $ is a
neighborhood of $0$ in $G,$ so $H$ is an open subgroup of $G$. At this point
to deduce that the subgroup $2H$ is also open, a little digression is
necessary.

The endomorphism $\theta :z\rightarrow 2z$ being continuous, its image $2G$
is a compact subgroup. Since $G$ is monothetic the index of the subgroup $2G$
in $G$ is $1$ or $2$ : $G=(2G)\cup (\alpha +2G)$ the two classes being
compact and equal or disjoint. Hence $2G$ is open. The subgroup $H$ being
open there is a smallest positive integer $j$ such that $j\alpha \in H$ and
the sequence $(nj\alpha )_{n}$ is dense in $H,$ hence $H$ is also
monothetic. Since $\theta $ is also an endomorphism of $H$ the subgroup $2H$
is open in $H$ hence in $G.$

The index $[G:2H]$ of the subgroup $2H$ in $G$ is finite since $2H$ is open
and $G$ is compact, and $m(2H)=1/[G:2H].$ The set $A^{a}$ is invariant by $%
2H $ hence it is union of classes of $2H$ and we obtain $\mu \left\{ f\leq
a\right\} =m(A^{a})\geq 1/[G:2H].$

This is a contradiction since the subgroup $2H$ depends only on the set $E$
and not on $a,$ and, as was put at the beginning, $\mu \left\{ f\leq
a\right\} $ can be positive and arbitrarily small. Therefore $\mu
(E)=0.\medskip $

If instead of $f(T^{-n}x)<f(T^{n}x)$ we consider inequalities

\noindent $f(T^{-n}x)<f(T^{n+q}x)$ with $q$ fixed, the proof has to be
slightly modified. The set $E$ is replaced by 
\begin{equation*}
E_{q}=\cup _{N>0}\cap _{n>N}\left\{ f\circ T^{-n}<f\circ T^{n+q}\right\}
\end{equation*}
and the sets $E_{N}^{a}$ by 
\begin{equation*}
E_{N,q}^{a}=\cap _{n>N}(T^{n}\left\{ f\leq a\right\} \cup T^{-n-q}\left\{
f>a\right\} ).
\end{equation*}
The key equality

\noindent $\int_{G}\varphi (y)\int_{G}\psi ^{a}(z)(1-\psi
^{a}(2y-z))dm(z)dm(y)=0$ becomes%
\begin{equation*}
\int_{G}\varphi _{q}(y)\int_{G}\psi ^{a}(z)(1-\psi ^{a}(2y+q\alpha
-z))dm(z)dm(y)=0
\end{equation*}%
with $\varphi _{q}=\pi 1_{E_{q}}$, since in Furstenberg's identity we have
to replace $1_{\left\{ f\leq a\right\} }=w$ by $1_{\left\{ f\leq a\right\}
}\circ T^{q}$ hence $\psi ^{a}=\pi w$ by $\psi ^{a}\circ \tau _{q\alpha }.$
The analysis of this new equality remains almost the same: instead of $%
\left\{ \psi ^{a}>0\right\} \subset (-A^{a})+2y$ we get $\left\{ \psi
^{a}>0\right\} \subset (-A^{a})+2y+q\alpha $ and when we take differences of
transformations $z\rightarrow 2y+q\alpha -z,$ the term $q\alpha $
disappears. So Theorem 5 is completely proved.$\diamond \medskip $

\textbf{Remarks}. If the system is weakly mixing a shorter proof of Theorem
5 is possible.

In theorems 4 and 5 the strict inequalities cannot be replaced by wide ones.
As an example consider the group $G=\mathbb{Z}/3\mathbb{Z}\times \mathbb{R}/%
\mathbb{Z}$ with $\alpha =$ $(1,\epsilon )$ where $\epsilon $ is irrational;
take $f(u,x)=1$ if $u=0,$ $0$ if $u=1,$ $-1$ if $u=2.$ Then it is clear that 
$B_{n}f=0$ on $\left\{ 1\right\} \times \mathbb{R}/\mathbb{Z}$ for every $%
n\geq 0.$

In this work the assumption of ergodicity was only a simplification. For a
non ergodic measure preserving dynamical system defined on a Lebesgue space,
the results of this paper are easily reformulated using the ergodic
decomposition of the invariant measure; in particular, without any change,
Theorem 5 still holds.$\medskip $

\textbf{Acknowledgement}. The authors thank the referee for her or his
careful reading of the paper and her or his pertinent suggestions.\bigskip

\textbf{References\medskip }

[1] Y. Derriennic, \textit{Ergodic theorem, reversibility and the filling
scheme}, Colloq. Math. 118 (2010), 599-608.\medskip

[2] H. Furstenberg, \textit{Ergodic behavior of diagonal measures and a
theorem of Szemer\'{e}di on arithmetic progressions}, Journal D'analyse Math%
\'{e}matique, Vol. 31(1977), 204-256.\medskip

[3] G. Hal\'{a}sz, \textit{Remarks on the remainder in Birkhoff's ergodic
theorem,} Acta Math. Acad. Sci. Hungar. 28 (1976), 389-395.\medskip

[4] M. Hochman, \textit{On one-sided and two-sided ergodic averages of
non-integrable functions}, Einstein Institute Note (2010).\medskip

[5] U. Krengel, Ergodic theorems, de Gruyter Studies in Math. n$%
{{}^\circ}%
6$ (1985).\medskip

[6] K. Petersen, Ergodic theory, Cambridge Studies in Math. n$%
{{}^\circ}%
2$ (1983).\medskip

[7] J. Wo\'{s}, \textit{Approximate convergence in ergodic theory}, Proc.
London Math. Soc. (3) 53(1986), 65-84.

\bigskip

Christophe.Cuny@univ-brest.fr

Yves.Derriennic@univ-brest.fr

Universit\'{e} de Brest

Math\'{e}matiques (LMBA-UMR CNRS 6205)

6 av. Le Gorgeu

29238 Brest F.

\end{document}